\documentclass{article}
\usepackage{graphicx, amsmath,amssymb,pxfonts}
\newtheorem{thm}{Theorem}
\newtheorem{lem}[thm]{Lemma}
\newtheorem{rem}[thm]{Remark}
\title{Efficient Discretization of  Stochastic Integrals}
\author{Masaaki Fukasawa \\
               Department of Mathematics, Osaka University 
}
\date{}
\begin{document}
\maketitle
\begin{abstract}
Sharp asymptotic lower bounds of the expected 
quadratic variation of 
discretization error in 
 stochastic integration are given. 
The theory relies on  inequalities for the kurtosis and skewness
 of a general random variable which are themselves seemingly new.
Asymptotically efficient schemes 
which attain the lower bounds are constructed explicitly.
The result is directly applicable to  practical hedging
problem in mathematical finance; it gives an  asymptotically optimal way to
 choose rebalancing dates and portofolios with respect to transaction costs.
The asymptotically 
efficient strategies in fact reflect the structure of transaction costs.
In particular a specific biased rebalancing scheme is
 shown to be superior to unbiased schemes if transaction costs  follow
 a convex model.
The problem is discussed also in terms of the exponential utility maximization.
\end{abstract}

\section{Introduction}
The stochastic integral $X \cdot Y_\sigma$ with respect to a 
semimartingale $Y$ and a stopping time $\sigma$
is by definition a limit of
$X^n \cdot Y_{\sigma}$ in probability, where $X^n$ is a sequence of
simple predictable processes with
$\sup_{t \in [0,\sigma]}|X^n_t - X_t| \to 0$
in probability as $n\to \infty$.
This convergence of $X^n \cdot Y$ is essential not only for the
theoretical construction
of the stochastic integral but also for  practical approximations
 in problems modeled with stochastic integrals.
The aim of this paper is to give a way to choose $X^n$ efficiently
in an asymptotic sense.
The main assumption of the paper is 
that $X$ is a continuous semimartingale.

Denote by $K$ the Radon-Nikodym derivative of the absolutely continuous
part of $\langle Y \rangle$ with respect to $\langle X \rangle$,
which always exists in light of the Lebesgue decomposition theorem.
Fukasawa~\cite{SAFA} showed that
\begin{equation*}
\liminf_{n \to \infty}
\mathbb{E}[N[X^n]_\sigma] \mathbb{E}[\langle Z[X^n] \rangle_\sigma] 
\geq \frac{1}{6} \mathbb{E}[\sqrt{K} \cdot \langle X \rangle_\sigma],
\end{equation*}
where $N[X^n]_\sigma$ is the number of the jumps of a given 
simple predictable process 
$X^n$ up to $\sigma$ and
$Z[X^n]:= (X-X^n) \cdot Y$ is the associated approximation error.
If $Y$ is a local martingale, then 
$\mathbb{E}[|Z[X^n]_{\sigma}|^2] =\mathbb{E}[\langle Z[X^n]
\rangle_\sigma]$ under a reasonable assumption,
and so the above inequality gives 
an asymptotic  lower bound of the mean squared error of discretization.
Notice that the bound does not depend on $X^n$.
The inequality is sharp in that the lower bound is attained by
\begin{equation} \label{eff1}
\begin{split}
 &X^n_t := X_{\tau^n_j}, \ \ t \in (\tau^n_j, \tau^n_{j+1}], \ \
 j=0,1,\dots,\\
&\tau^n_0 := 0, \ \ \tau^n_{j+1} := \inf\{t > \tau^n_j ; |X_t
 -X_{\tau^n_j}| = \epsilon_n K_{\tau^n_j}^{-1/4}\}, \ \ \epsilon_n
 \downarrow 0
\end{split}
\end{equation}
under a reasonable condition.
We call such a sequence $X^n$ that attains the lower bound an
asymptotically efficient scheme.
The above result is extended and proved under a less restrictive condition in
this paper as a particular case.

To obtain a precise approximation to $X \cdot Y$, one has to take $X^n$
as close to $X$ as possible. In practical contexts it may be inevitably
accompanied by various kinds of cost, especially if $X$ is not of
finite variation. The number of jumps $N[X^n]_\sigma$ is
interpreted as one of them.
In the context of mathematical finance for example,
$X$ and $Y$ stand for a portfolio strategy and an asset price
 process respectively.
Then $Z[X^n]$ represents the replication error associated to a discrete 
rebalancing strategy $X^n$.
A continuous rebalancing is impossible in practice and $N[X^n]_\sigma$
corresponds to the number of trading, a measure on trader's effort.
The scheme (\ref{eff1}) defines an asymptotically
efficient discrete strategy
which asymptotically minimizes the mean squared error relative to 
the specific cost function $\mathbb{E}[N[X^n]_\sigma]$.

The sequence  $\mathbb{E}[N[X^n]_\sigma]$ is however just one of
measures on costs.
Again for example in the financial context, 
the cumulative transaction cost associated to $X^n$ is often modeled as
\begin{equation*}
\kappa \sum_{0 < t \leq \sigma} Y_t | \Delta X^n_t|
\end{equation*}
with a constant $\kappa > 0$.
This is the so-called linear or proportional transaction cost
model. 
More generally one may consider 
as a cost or penalty,
\begin{equation} \label{cost}
C[S,\beta; X^n] := \sum_{0 < t \leq \sigma}S_tK_t|\Delta X^n_t|^\beta
1_{\{|\Delta X^n| > 0\}}
\end{equation}
with a nonnegative predictable process $S$ and
a constant  $\beta \geq 0$. 
Notice  that $C[1/K,0;X^n]_\sigma$ and $C[Y/K,1;X^n]_{\sigma}$
represent the number of rebalancing and  the cumulative 
linear transaction cost
respectively.
If $\beta \in (0,1)$ or $\beta >1 $, the cost is concave or convex
respectively in the amount of transaction.
Beyond these interpretations in the financial context, we treat
the general form of $C[S,\beta;X^n]_\sigma$ as a penalty 
against taking $X^n$ too close to $X$.
Then a natural problem would be to minimize 
$\mathbb{E}[\langle Z[X^n] \rangle_\sigma]$ 
relative to the expected cost $\mathbb{E}[C[S,\beta; X^n]_\sigma]$
in the asymptotic  situation that
$\sup_{t \in [0,\sigma]}|X^n_t - X_t| \to 0.$
Fukasawa~\cite{RFE}(in Japanese) proposed this framework and proved that
for all $\beta \in [0,2)$,
\begin{equation} \label{cx}
\liminf_{n \to \infty}
|\mathbb{E}[C[S,\beta;X^n]_\sigma]|^{2/(2-\beta)} \mathbb{E}[\langle Z[X^n] \rangle_\sigma] 
\geq \frac{1}{6} |\mathbb{E}[ (S^{2/(4-\beta)}K) \cdot \langle X\rangle_\sigma]|^{(4-\beta)/(2-\beta)}
\end{equation}
if $X^n$ is of the form
$X^n_t = X_{\tau^n_j}$ for any
 $t \in (\tau^n_j,\tau^n_{j+1}]$ 
with 
 an increasing sequence of stopping times $\tau^n = \{\tau^n_j\}$ with
$\tau^n_0 = 0$ and
$\sup_{j \geq 0}|\tau^n_{j+1} \wedge \sigma - \tau^n_j \wedge \sigma|  \to 0$ as 
$n\to \infty$.
The lower bound is sharp in that it is attained by
\begin{equation} \label{eff}
\begin{split}
&X^n_t := X_{\tau^n_j}, \ \ t \in (\tau^n_j, \tau^n_{j+1}], \ \
 j=0,1,\dots,\\
&\tau^n_0=0, \ \ \tau^n_{j+1} =
\inf\left\{t > \tau^n_j; |X_t - X_{\tau^n_j}| \geq 
\epsilon_n S_{\tau^n_j}^{1/(4-\beta)}\right\}, \ \ \epsilon_n \downarrow
 0
\end{split}
\end{equation}
under a reasonable condition.
The proof is given in this paper as well under a less restrictive condition.
This result does not give a complete answer to our problem in that the
lower bound is for a restricted class of $X^n$ as 
$X^n_t = X_{\tau^n_j}$ for $t \in (\tau^n_j,\tau^n_{j+1}]$ 
with some $\{\tau^n_j\}$.
We call such $X^n$ an unbiased scheme. 
Intuitively, taking $X^n$ in the unbiased manner is natural and
necessary to have a good approximation to $X \cdot Y$.
In fact in the case $\beta = 0$ and
$C[S,\beta,X;X^n] = N[X^n]_\sigma$,
as stated first, the unbiased scheme $X^n$ defined by (\ref{eff1}) is 
asymptotically efficient.
The main result of this paper shows that the discretization scheme
(\ref{eff}) is actually asymptotically 
efficient if $\beta \in [0,1]$, however not
so if $\beta \in (1,2)$. In the latter case, surprisingly,
the lower bound is reduced
to one third and asymptotically attained by a sequence of biased schemes.

In Section~2, we give a general result on the centered moments of a
random variable, which seems new and important itself and 
plays  an essential role to derive lower 
bounds of discretization error in the stochastic integration.
In Section~3, we give a sharp lower bound for unbiased schemes, which is
a slight extension of the result of Fukasawa~\cite{RFE}(in Japanese).
In Section~4, we give sharp lower bounds for possibly biased schemes and
construct explicit schemes which asymptotically attain the bounds. 
 In Section~5, 
we show that an asymptotically 
efficient scheme  is a maximizer of a scaling limit of the
exponential utility  in the financial context of discrete hedging.

We conclude this section by mentioning related studies in the
literature.
Rootz\'en~\cite{Rootzen} studied the discretization error of stochastic
integrals with the equidistant partition $\tau^n_j = j/n$ and 
proved that the discretization error of
a stochastic integral converges in law to a time-changed Brownian motion 
with rate $n^{-1/2}$ as $n \to \infty$.
An extension to discontinuous semimartingales was given by Tankov and
Voltchkova~\cite{TV} in the equidistant case.
Fukasawa~\cite{AAP} gave an extension to another direction that admits a general sequence of locally
homogeneous  stochastic partitions
and gave several sharp
lower bounds of the asymptotic conditional variance of the
discretization error.
Hayashi and Mykland~\cite{HM2005} revisited Rootz\'en's problem 
in terms of the discrete hedging in mathematical finance.
Motivated by this financial application, the mean squared error
was studied by Gobet and Temam~\cite{GT01}, Geiss and Geiss~\cite{GG06},
Geiss and Toivola~\cite{GT09} under the Black-Scholes model.
Among others, Geiss and Geiss~\cite{GG06} showed that the use of 
stochastic partitions does not improve the convergence rate. In a sense
our result refines this observation under a general framework.
Our problem is also related to Leland's strategy for hedging under
transaction costs. See Leland~\cite{Leland}, Denis and
Kabanov~\cite{DK},
Fukasawa~\cite{RAFE}.
The difference is that we are looking for an efficient 
discrete hedging strategy which does not require a surcharge, 
while Leland's strategy does it to absorb transaction costs.
In a statistical framework, Genon-Catalot and Jacod~\cite{GJ}
studied an optimality problem for a class of random sampling schemes,
which is smaller than our class.
Finally remark that the use of hitting times such as (\ref{eff})
has another advantage in terms of almost sure convergence.
See Karandikar~\cite{Karandikar1995}.

\section{Kurtosis-skewness inequalities}
Here we study the centered moments of a general random variable.
The reason why we need such a general framework is that 
in our problem of discretization, we encounter the moments of
a martingale evaluated at a  stopping time, 
which can follow any distribution with mean $0$ in light of Skorokhod
stopping problem.
The notation in this section  is independent of that in other sections. 
We say a random variable $X$ is Bernoulli if the support of $X$ consists
of two points.
We say $X$ is symmetrically Bernoulli if $X$ is
Bernoulli and its skewness is $0$, that is, 
$\mathbb{E}[(X-\mathbb{E}[X])^3] = 0$.
For any random variable $X$ with $\mathbb{E}[X] = 0$,
$\mathbb{E}[X^2] > 0$ and $\mathbb{E}[X^4] < \infty$,
it holds that
\begin{equation} \label{Pearson}
\frac{\mathbb{E}[X^4]}{|\mathbb{E}[X^2]|^2} - 
\frac{|\mathbb{E}[X^3]|^2}{|\mathbb{E}[X^2]|^3} \geq 1.
\end{equation}
This is often called Pearson's inequality and easily shown as follows:
\begin{equation*}
|\mathbb{E}[X^3]|^2 = 
|\mathbb{E}[X(X^2-\mathbb{E}[X^2])]|^2 \leq
\mathbb{E}[X^2] (\mathbb{E}[X^4]- |\mathbb{E}[X^2]|^2).
\end{equation*}
From this proof it is clear that
the equality is attained only if $X$ is Bernoulli.
Conversely if $X$ is Bernoulli, then we get the equality by a
straightforward calculation.
Pearson's inequality was used by Fukasawa~\cite{AAP}\cite{SAFA} to obtain lower
bounds of discretization error of stochastic integrals.
This is however not sufficient for our current purpose.
Fukasawa~\cite{AAP} proved another inequality which looks similar to 
but independent of (\ref{Pearson}):
\begin{equation} \label{Fukasawa}
\frac{\mathbb{E}[X^4]}{|\mathbb{E}[X^2]|^2} - \frac{3}{4}
\frac{|\mathbb{E}[X^3]|^2}{|\mathbb{E}[X^2]|^3} \geq 
 \frac{\mathbb{E}[X^2]}
{|\mathbb{E}[|X|]|^2}.
\end{equation}
The equality is attained if and only if $X$ is Bernoulli.
The proof is lengthy and unexpectedly different from that for Pearson's
inequality. See Appendix~B of Fukasawa~\cite{AAP}. 
From these inequalities we obtain the following lemmas.
\begin{lem}\label{KS1}
Let $\beta \in [0,1)$.
For any random variable $X$ with $\mathbb{E}[X] = 0$,
$\mathbb{E}[X^2] > 0$ and $\mathbb{E}[X^4] < \infty$,
\begin{equation} \label{ks1}
\frac{\mathbb{E}[X^4]}{|\mathbb{E}[X^2]|^2} - \frac{3}{4}
\frac{|\mathbb{E}[X^3]|^2}{|\mathbb{E}[X^2]|^3} \geq 
 \frac{|\mathbb{E}[X^2]|^{\beta/(2-\beta)}}
{|\mathbb{E}[|X|^{\beta}]|^{2/(2-\beta)}}.
\end{equation}
The equality is attained if and only if $X$ is symmetrically
 Bernoulli.
\end{lem}
{\it Proof: } By H$\ddot{\text{o}}$lder's inequality, we have
\begin{equation*}
\mathbb{E}[|X|] \leq |\mathbb{E}[X^2]|^{(1-\beta)/(2-\beta)}
|\mathbb{E}[|X|^\beta]|^{1/(2-\beta)},
\end{equation*}
or equivalently,
\begin{equation*}
\frac{\mathbb{E}[X^2]}{|\mathbb{E}[|X|]|^2}
\geq \frac{|\mathbb{E}[X^2]|^{\beta/(2-\beta)}}
{|\mathbb{E}[|X|^{\beta}]|^{2/(2-\beta)}}.
\end{equation*}
The result then follows from (\ref{Fukasawa}). \hfill////
\begin{lem}
Let $\beta \in [0,2)$ and $\alpha \in [0,1]$.
For any random variable $X$ with $\mathbb{E}[X] = 0$,
$\mathbb{E}[X^2] > 0$ and $\mathbb{E}[X^4] < \infty$,
\begin{equation} \label{ks20}
\frac{\mathbb{E}[X^4]}{|\mathbb{E}[X^2]|^2} - \alpha
\frac{|\mathbb{E}[X^3]|^2}{|\mathbb{E}[X^2]|^3} - (1-\alpha)
 \frac{|\mathbb{E}[X^2]|^{\beta/(2-\beta)}}
{|\mathbb{E}[|X|^{\beta}]|^{2/(2-\beta)}}  \geq \alpha.
\end{equation}
The equality is attained if and only if $X$ is symmetrically
 Bernoulli.
\end{lem}
{\it Proof: } By H$\ddot{\text{o}}$lder's inequality, we have
\begin{equation*}
\mathbb{E}[X^2] \leq |\mathbb{E}[|X|^\beta]|^{2/(4-\beta)}
|\mathbb{E}[X^4]|^{(2-\beta)/(4-\beta)},
\end{equation*}
or equivalently,
\begin{equation*}
 \frac{\mathbb{E}[X^4]}
{|\mathbb{E}[X^2]|^2}
\geq \frac{|\mathbb{E}[X^2]|^{\beta/(2-\beta)}}
{|\mathbb{E}[|X|^{\beta}]|^{2/(2-\beta)}}.
\end{equation*}
Therefore,
\begin{equation*}
\frac{\mathbb{E}[X^4]}{|\mathbb{E}[X^2]|^2} - \alpha
\frac{|\mathbb{E}[X^3]|^2}{|\mathbb{E}[X^2]|^3} - (1-\alpha)
 \frac{|\mathbb{E}[X^2]|^{\beta/(2-\beta)}}
{|\mathbb{E}[|X|^{\beta}]|^{2/(2-\beta)}}  \geq 
\alpha 
\left\{
\frac{\mathbb{E}[X^4]}{|\mathbb{E}[X^2]|^2} - 
\frac{|\mathbb{E}[X^3]|^2}{|\mathbb{E}[X^2]|^3}
\right\}.
\end{equation*}
The result then follows from (\ref{Pearson}). \hfill////
\begin{lem}\label{KS2}
Let $\beta \in [0,2)$ and $\alpha \in (0,1]$.
For any random variable $X$ with $\mathbb{E}[X] = 0$,
$\mathbb{E}[X^2] > 0$ and $\mathbb{E}[X^4] < \infty$,
\begin{equation} \label{ks2}
 \frac{|\mathbb{E}[|X|^{\beta}]|^{2/(2-\beta)}}
{|\mathbb{E}[X^2]|^{\beta/(2-\beta)}}
\left\{\frac{\mathbb{E}[X^4]}{|\mathbb{E}[X^2]|^2} - \alpha
\frac{|\mathbb{E}[X^3]|^2}{|\mathbb{E}[X^2]|^3} \right\} > 1-\alpha.
\end{equation}
Moreover if $X$ is Bernoulli,  then 
\begin{equation} \label{ks2s}
\frac{|\mathbb{E}[|X|^{\beta}]|^{2/(2-\beta)}}
{|\mathbb{E}[X^2]|^{\beta/(2-\beta)}}
\left\{\frac{\mathbb{E}[X^4]}{|\mathbb{E}[X^2]|^2} - \alpha
\frac{|\mathbb{E}[X^3]|^2}{|\mathbb{E}[X^2]|^3} \right\} = 
F_{\alpha\beta}\left(\frac{|\mathbb{E}[X^3]|^2}{|\mathbb{E}[X^2]|^3}\right),
\end{equation}
where $F_{\alpha\beta}$ is a continuous function
 with $F_{\alpha\beta}(0) = 1$.
If $\beta \in (1,2)$, then 
$F_{\alpha \beta}(\infty) = 1-\alpha$.
\end{lem}
{\it Proof: }
The inequality (\ref{ks2}) is apparent from (\ref{Pearson}) and 
(\ref{ks20}).
Let $X$ be Bernoulli.
We suppose $\mathbb{E}[X^2] = 1$ without loss of generality.
Then the support of $X$ is of the form $\{e^x,-e^{-x}\}$ 
and $\mathbb{P}[X = e^x] = 1/(1+e^{2x})$
with $x \in \mathbb{R}$.
By a straightforward calculation, we get
$\mathbb{E}[X^3] = 2 \sinh(x)$ and
\begin{equation*}
 \frac{|\mathbb{E}[|X|^{\beta}]|^{2/(2-\beta)}}
{|\mathbb{E}[X^2]|^{\beta/(2-\beta)}}
\left\{\frac{\mathbb{E}[X^4]}{|\mathbb{E}[X^2]|^2} - \alpha
\frac{|\mathbb{E}[X^3]|^2}{|\mathbb{E}[X^2]|^3} \right\}
 = \frac{4 \alpha -3 + 4(1-\alpha)|\cosh(x)|^2}
{|\cosh(x)|^{2/(2-\beta)} |\cosh((\beta-1)x)|^{-2/(2-\beta)} }.
\end{equation*}
Putting
\begin{equation} \label{gdef}
g(x) = \cosh((\beta-1)x)|\cosh(x)|^{1-\beta},
\end{equation}
the right hand side is given by
\begin{equation}\label{rhs}
\frac{4\alpha-3}{g(x)^{-2/(2-\beta)} |\cosh(x)|^2}
+ \frac{4(1-\alpha)}{g(x)^{-2/(2-\beta)}}.
\end{equation}
Notice that $g(0) = 1$ and 
$g(x)^{-2/(2-\beta)}$ converges to $4$ as $|x| \to \infty$
for $\beta \in (1,2)$.
\hfill////
\begin{rem} \upshape
 Let $g$ be defined by (\ref{gdef}).
Since
\begin{equation*}
g^\prime(x) = (\beta-1)\sinh((\beta-2)x) |\cosh(x)|^{-\beta}, 
\ \ g^{\prime \prime}(0) = (1-\beta)(2-\beta),
\end{equation*}
for $\beta \neq 1$,  $g^\prime(x) = 0$ if and only if $x=0$.
Further if $\beta \in [0,1)$ or 
$\beta \in (1,2)$, respectively, 
the minimum or maximum of $g$ is attained at $x =0$.
Therefore if $\alpha \geq 3/4$ and $\beta \in (1,2)$, 
the function defined by (\ref{rhs}) is decreasing in $|x|$ 
and converges to $1-\alpha$ as $|x| \to \infty$.
However in the following
 sections, we use Lemma~\ref{KS2} with $\alpha = 2/3$,
 where the function is  not necessarily monotone in $|x|$.
\end{rem}

\section{Efficiency for unbiased Riemann sums}
Here we recall  the problem with a rigorous formulation and
give a slight improvement of the result of Fukasawa~\cite{RFE}.
Let $X$ and $Y$ be semimartingales defined on a filtered probability
space $(\Omega, \mathcal{F}, \mathbb{P}, \{\mathcal{F}_t\})$ which satisfies
the usual conditions.
We assume that 
that there exist a continuous local martingale $M$ and
a locally bounded adapted process $H$ such that
\begin{equation*}
X = H \cdot \langle M \rangle + M.
\end{equation*}
Denote by  $\mathcal{T}$ the set of the  
increasing sequences of stopping times $\tau = \{\tau_j\}$ with
$0 = \tau_0 < \tau_1 < \cdots$ and 
$\lim_{j \to \infty} \tau_j = \infty$ a.s.. 
Given $\tau = \{\tau_j\}\in \mathcal{T}$,
define a simple predictable process 
$X[\tau]$ as $X[\tau]_t = X_{\tau_j}$ for
$t \in (\tau_j, \tau_{j+1}]$.
Conversely, for a given simple predictable process $\hat{X}$, 
define $\tau[\hat{X}] \in \mathcal{T}$ as the sequence of the 
jump times of $\hat{X}$.
By definition we have
\begin{equation*}
\begin{split}
&Z[X[\tau]]_t = \int_0^t X_s \mathrm{d}Y_s - \sum_{j=0}^{\infty}
X_{\tau_j}(Y_{\tau_{j+1} \wedge t} -Y_{\tau_j\wedge t} ),\\
&Z[\hat{X}]_t = \int_0^t X_s \mathrm{d}Y_s - \sum_{j=0}^{\infty}
\hat{X}_{\tau[\hat{X}]_j + }(Y_{\tau[\hat{X}]_{j+1} \wedge t} -
Y_{\tau[\hat{X}]_j\wedge t} )
\end{split}
\end{equation*}
for $t \geq 0$. Our aim is to minimize
$\mathbb{E}[\langle Z[X^n] \rangle_{\sigma}]$ asymptotically when
\begin{equation} \label{sup}
\sup_{t \in [0,\sigma]}|X^n_t - X_t| \to 0
\end{equation}
in probability as $n \to \infty$. 
Denote by $K$ the Radon-Nikodym derivative of the absolutely continuous
part of the predictable quadratic variation 
$\langle Y \rangle$ with respect to $\langle X \rangle$,
which always exists in light of the Lebesgue decomposition theorem.
We consider the cost $C[S,\beta; \hat{X}]$ defined by (\ref{cost}) for a
given simple predictable process $\hat{X}$.
We assume that $K$ and $S$ are positive, 
continuous and moreover, constant on any random
interval where $\langle X \rangle$ is constant.
By the last assumption, we have
\begin{equation}\label{Xconti}
K = \tilde{K}_{\langle X \rangle}, \ \ \tilde{K} = K_{F}, \ \ 
S = \tilde{S}_{\langle X \rangle}, \ \ \tilde{S} = S_{F},
\end{equation}
where $F_s = \inf\{t \geq 0;\langle X \rangle_t > s\}$;
see Karatzas and Shreve~\cite{KS}, 3.4.5. 

Now we define a class of unbiased schemes in which at first we consider the
efficiency or optimality of discretization.
Denote by $\mathcal{T}_u(S,\beta,\sigma)$ the set of the sequences
of simple predictable processes $X^n$ of the form $X^n = X[\tau^n]$,
$\tau^n = \{\tau^n_j\}\in \mathcal{T}$, such that 
there exists a sequence of stopping times $\sigma^m$ with 
$\sigma^m \to \sigma$ as $m \to \infty$,
\begin{enumerate}
\item for each $m$, (\ref{sup}) holds with $\sigma^m$ instead of $\sigma$.
\item  for each $m$,
\begin{equation*}
\mathbb{E}[C[S,\beta;X^n]_{\sigma^m}]^{2/(2-\beta)}
\langle Z[X^n]\rangle_{\sigma^m}
\end{equation*}
is uniformly integrable in $n$.
\end{enumerate}
\begin{rem} \upshape
The uniform integrability condition for $\mathcal{T}_u(S,\beta,\sigma)$
is usually easy to check.
It is for example satisfied when considering  the sequence of the
equidistant partitions $\tau^n_j = j/n$ if
$ \mathrm{d}\langle X \rangle_t$ has a 
locally bounded Radon-Nikodym derivative with respect to $\mathrm{d}t$.
The exponent $2/(2-\beta)$ is actually chosen so that
$\mathbb{E}[C[S,\beta;X^n]_{\sigma^m}]^{2/(2-\beta)}  \propto n $
asymptotically in the equidistant case since $n^{-1}$ is the optimal
convergence rate of $\langle Z[X^n]\rangle_{\sigma^m}$ for the case.
All reasonable $X^n$ should enjoy this property of rate-efficiency.
Note that by the Dunford-Petis theorem, the uniform integrability 
is equivalent to
the relative compactness in the 
$\sigma(L^1,L^\infty)$ topology.
By the Eberlein-Smulian theorem,
it is further equivalent to 
the relative sequential compactness in the same topology.
\end{rem}

\begin{thm} \label{R1}
Let $\beta \in [0,2)$.
The inequality (\ref{cx}) holds for 
 all $\{X^n\} \in \mathcal{T}_u(S,\beta,\sigma)$.
\end{thm}
For the proof, we start with a lemma.
\begin{lem} \label{supsup}
Let $X^n$ be a sequence of simple predictable processes.
Then  (\ref{sup}) implies that
\begin{equation}\label{brasup}
\sup_{j \geq 0} |\langle X \rangle_{\tau^n_{j+1} \wedge \sigma}
 - \langle X \rangle_{\tau^n_j \wedge
 \sigma}| \to 0
\end{equation}
in probability as $n \to \infty$ with $\tau^n = \tau[X^n]$.
Conversely if (\ref{brasup}) holds  for a sequence 
$\tau^n \in  \mathcal{T}$, then (\ref{sup}) 
holds with $X^n =
 X[\tau^n]$.
\end{lem}
{\it Proof: }
For any subsequence of $n$, there exists a further subsequence $n_k$
such that (\ref{sup}) holds a.s. with $n=n_k$ as $k \to \infty$.  
It suffices then to show that (\ref{brasup}) holds a.s. with this
subsequence.
Let $\Omega^\ast$ be a subset of $\Omega$ such that
for any $\omega \in \Omega^\ast$, (\ref{brasup}) does not hold 
with $n = n_k$, $k\to \infty$.
Then, for $\omega \in \Omega^\ast$,
there exist $\epsilon(\omega)>0$ and a sequence of intervals 
$I_m(\omega) = [a_m(\omega),b_m(\omega)]$ such that
for each $m$, there exists $n=n_k$ such that
$I_m(\omega) = [\tau^n_{j}(\omega), \tau^n_{j+1}(\omega)]$ and
\begin{equation*}
\inf_m |\langle X \rangle_{b_m}(\omega) - \langle X
 \rangle_{a_m}(\omega)|
\geq \epsilon(\omega).
\end{equation*}
Since $(a_m(\omega), b_m(\omega))$ is a sequence in the compact set
$[0,\sigma(\omega)] \times [ 0, \sigma(\omega)]$,
it has an accumulating point $[a_\ast(\omega),b_\ast(\omega)]$ with
\begin{equation*}
 |\langle X \rangle_{b_\ast}(\omega) - \langle X
 \rangle_{a_\ast}(\omega)|
\geq \epsilon(\omega).
\end{equation*}
With probability one,
$\langle X \rangle$ is continuous, so we may suppose that
$a_\ast(\omega) < b^\ast(\omega)$ without loss of generality.
Again with probability one,
if $X$ is constant on an interval, then
 $\langle X \rangle$ is constant on the interval.
So we may suppose that $X(\omega)$ is not constant on 
$[a_\ast(\omega), b_\ast(\omega)]$ without loss of generality.
On the other hand, there exists a subsequence $X^m(\omega)$ of 
$X^{n_k}(\omega)$ 
such that $X^m(\omega)$ is constant on a non-empty interval of
$[a_\ast(\omega), b_\ast(\omega)]$.
Recalling the way that the subsequence was chosen,
we conclude that $\mathbb{P}[\Omega^\ast] = 0$.
\hfill////

\vspace*{1cm}
\noindent
{\it Proof of Theorem~\ref{R1}: }
Put $\tau^n = \tau[X^n]$.
By the usual localization argument, we may and do suppose
without loss of generality that
$X,  \langle X \rangle, K, 1/K, S$ and $H$
are bounded up to $\sigma$,  that
(\ref{brasup}) holds, and that
$|\mathbb{E}[C[S,\beta;X^n]_{\sigma}]|^{2/(2-\beta)}
\langle Z[X^n]\rangle_{\sigma}$
is uniformly integrable in $n$.
Define 
$K[\tau]$
as $K[\tau]_t = K_{\tau_j}$ for 
$t \in [\tau_j, \tau_{j+1})$ for $\tau \in \mathcal{T}$.
Let
\begin{equation*}
\epsilon^n = \sup_{0 \leq s \leq \sigma}|K_s -
K[\tau^n]_s|.
\end{equation*}
By Lemma~\ref{supsup} and (\ref{Xconti}), we have that
$\epsilon_n$ is bounded and converges to $0$ in probability 
as $n\to \infty$.
By It$\hat{\text{o}}$'s formula,
\begin{equation} \label{ito}
\begin{split}
\langle Z[X^n] \rangle_t =& 
\int_0^t (X_s - X^n_s)^2 \mathrm{d} \langle Y \rangle_s \\
\geq &
\int_0^t (X_s - X^n_s)^2 K_s \mathrm{d} \langle X \rangle_s \\
=&
\int_0^t (X_s - X^n_s)^2 K[\tau^n]_s \mathrm{d} \langle X \rangle_s 
+
\int_0^t (X_s - X^n_s)^2 (K_s - K[\tau^n]_s
)\mathrm{d} \langle X \rangle_s  \\
=&
\frac{1}{6}\sum_{j=0}^{\infty}
K_{\tau^n_j} (X_{\tau^n_{j+1}\wedge t} - 
X_{\tau^n_{j}\wedge t} )^4 - 
\frac{2}{3} \int_0^t K[\tau^n]_s (X_s-X^n_s)^3 \mathrm{d}X_s
\\ &+ 
\int_0^t (X_s - X^n_s)^2 (K_s - K[\tau^n]_s
)\mathrm{d} \langle X \rangle_s.
\end{split}
\end{equation}

Now we show that
\begin{equation*}
\lim_{n \to \infty} \mathbb{E}\left[
|\mathbb{E}[C[S,\beta;X^n]_{\sigma}]|^{2/(2-\beta)}
\int_0^\sigma (X_s - X^n_s)^2 (K_s - K[\tau^n]_s)
\mathrm{d} \langle X \rangle_s  \right] = 0.
\end{equation*}
Put
\begin{equation*}
\begin{split}
V^n =& |\mathbb{E}[C[S,\beta;X^n]_{\sigma}]|^{2/(2-\beta)}
\langle Z[X^n] \rangle_{\sigma} \\ = &
|\mathbb{E}[C[S,\beta;X^n]_{\sigma}]|^{2/(2-\beta)}
\int_0^\sigma (X_s - X^n_s)^2 \mathrm{d} \langle Y \rangle_s.
\end{split}
\end{equation*}
Since $1/K$ is bounded by a constant, say, $A > 0$ and
$ K_s \mathrm{d}\langle X \rangle_s \leq \mathrm{d}\langle Y \rangle_s$,
we have
\begin{equation*}
|\mathbb{E}[C[S,\beta;X^n]_{\sigma}]|^{2/(2-\beta)}
\int_0^\sigma (X_s - X^n_s)^2 |K_s - K[\tau^n]_s|
\mathrm{d} \langle X \rangle_s 
\leq A  \epsilon^n V^n \to 0
\end{equation*}
in probability.
Since $\epsilon^n$ is bounded and $V^n$ is uniformly integrable, 
$\epsilon^n V^n$ is uniformly integrable as well and so, 
we obtain that $\mathbb{E}[\epsilon^n V^n] \to 0$.

Similarly, we can show that
\begin{equation*}
\begin{split}
&|\mathbb{E}[C[S,\beta;X^n]_{\sigma}]|^{2/(2-\beta)}
\mathbb{E}\left[
\int_0^\sigma K[\tau^n]_s(X_s-X^n_s)^3 \mathrm{d}X_s
\right]
\\&= 
|\mathbb{E}[C[S,\beta;X^n]_{\sigma}]|^{2/(2-\beta)}
\mathbb{E}\left[
\int_0^\sigma K[\tau^n]_s(X_s-X^n_s)^3 H_s \mathrm{d}\langle X
\rangle_s
\right]
\to 0
\end{split}
\end{equation*}
by using the continuity of $X$ instead of $K$.
So far we have obtained
\begin{equation*}
\begin{split}
&\liminf_{n\to \infty} 
|\mathbb{E}[C[S,\beta;X^n]_{\sigma}]|^{2/(2-\beta)}
\mathbb{E}[\langle Z[X^n] \rangle_\sigma] 
\\& \geq \liminf_{n \to \infty}
\frac{1}{6} 
|\mathbb{E}[C[S,\beta;X^n]_{\sigma}]|^{2/(2-\beta)}
\mathbb{E}\left[
\sum_{j=0}^{\infty}
K_{\tau^n_j} (X_{\tau^n_{j+1}\wedge \sigma} - 
X_{\tau^n_{j}\wedge \sigma} )^4 
\right].
\end{split}
\end{equation*}
On the other hand, 
by H$\ddot{\text{o}}$lder's inequality,
\begin{equation*}
\begin{split}
&\mathbb{E}
\left[
\sum_{j \geq 1, \tau^n_j \leq \sigma}
|S_{\tau^n_{j} }|^{2/(4-\beta)} |K_{\tau^n_{j-1}}|^{1/p}
|K_{\tau^n_{j} }|^{1/q}
(X_{\tau^n_{j} } - 
X_{\tau^n_{j-1}} )^2
\right] \\ &\leq 
\left|
\mathbb{E}
\left[
\sum_{j=0}^{\infty} 
K_{\tau^n_j}
 (X_{\tau^n_{j+1}\wedge \sigma} - 
X_{\tau^n_{j}\wedge \sigma} )^4 
\right]
\right|^{1/p}
\left|
\mathbb{E}
\left[
\sum_{0 < t \leq \sigma}
S_t K_t |\Delta X^n_t |^\beta 1_{\{|\Delta X^n_t| > 0\}}
\right]
\right|^{1/q}
\\ & =  
\left|
\mathbb{E}
\left[
\sum_{j=0}^{\infty} 
K_{\tau^n_j}
 (X_{\tau^n_{j+1}\wedge \sigma} - 
X_{\tau^n_{j}\wedge \sigma} )^4 
\right]
\right|^{1/p}
\left|
\mathbb{E}
\left[
C[S,\beta,X^n]_{\sigma}\right]
\right|^{1/q}
\end{split}
\end{equation*}
where $p = (4-\beta)/(2-\beta)$ and 
$q = p/(p-1) = (4-\beta)/2$.
The left hand side converges to 
$\mathbb{E}[(S^{2/(4-\beta)} K)\cdot \langle X
\rangle_{\sigma}]$.
\hfill////

\begin{thm} \label{R2}
Suppose that $\langle Y \rangle = K \cdot \langle X \rangle$.
Let $\hat{S}$ be a positive continuous adapted process which is constant
 on any random interval where $\langle X \rangle$ is constant.
Let $\epsilon_n$ be a positive sequence with $\epsilon_n \to 0$ 
as $n\to \infty$.
Define $X^n$ as
\begin{equation} \label{hitting} 
\begin{split}
&X^n_t := X_{\tau^n_j}, \ \ t \in [\tau^n_j, \tau^n_{j+1}), \ \
 j=0,1,\dots,\\
&\tau^n_0 := 0, \ \ \tau^n_{j+1} := \inf\{t > \tau^n_j ; |X_t
 -X_{\tau^n_j}| = \epsilon_n \hat{S}_{\tau^n_j}\}.
\end{split}
\end{equation}
Then $\{X^n\} \in  \mathcal{T}_u(S,\beta,\sigma)$ for any $\beta \in [0,2)$.  
Moreover if 
$X, \langle X \rangle,  H,  K,   1/K,  S, 1/S,  \hat{S}$ 
and $ 1/\hat{S}$
are bounded up to $\sigma$, then we have that for any $\beta \in [0,2)$,
\begin{equation*}
\sum_{j=0}^{\infty}
|X_{\tau^n_{j+1} \wedge \sigma} - X_{\tau^n_{j} \wedge \sigma}|^2
, \ \ 
\frac{C[S,\beta;X^n]_{\sigma}}
{\mathbb{E}[C[S,\beta;X^n]_{\sigma}]}
\end{equation*}
are uniformly integrable in $n$, and
\begin{equation*}
\begin{split}
&\lim_{n \to \infty}
 \epsilon_n^{2-\beta}
\mathbb{E}[C[S,\beta;X^n]_{\sigma}]|^{2/(2-\beta)}  = 
\mathbb{E}[(S\hat{S}^{\beta-2})\cdot \langle Y \rangle_\sigma]\\
& \lim_{n \to \infty}
 \epsilon_n^{-2} \mathbb{E}[\langle Z[X^n]\rangle_{\sigma}] = 
\frac{1}{6}\mathbb{E}[\hat{S}^2 \cdot \langle Y \rangle_\sigma].
\end{split}
\end{equation*}
In particular if $\hat{S} = S^{1/(4-\beta)}$, or equivalently,
 $X^n$ is defined by (\ref{eff}), then
\begin{equation*}
\lim_{n \to \infty}
|\mathbb{E}[C[S,\beta;X^n]_{\sigma}]|^{2/(2-\beta)}
\mathbb{E}[\langle Z[X^n]\rangle_{\sigma}]
= \frac{1}{6}
|\mathbb{E}[
(S^{2/(4-\beta)}K)  \cdot 
\langle X \rangle_{\sigma}]|^{(4-\beta)/(2-\beta)}.
\end{equation*}
\end{thm}
{\it Proof: }
By the usual localization argument, we may and do suppose
without loss of generality that
$X, \langle X \rangle,  H,  K,   1/K,  S, 1/S, \hat{S}$ and 
$1/\hat{S}$ 
are bounded up to $\sigma$.
Then, notice that
the uniformly integrability of
\begin{equation} \label{qv}
\sum_{j=0}^{\infty}
|X_{\tau^n_{j+1} \wedge \sigma} - X_{\tau^n_{j} \wedge \sigma}|^2
\end{equation}
follows from the decomposition
\begin{equation*}
\sum_{j=0}^{\infty}
|X_{\tau^n_{j+1} \wedge \sigma} - X_{\tau^n_{j} \wedge \sigma}|^2
 = \langle X \rangle_\sigma + 
2\int_0^\sigma (X_t-X^n_t)H_t\mathrm{d}\langle X \rangle_t
+ 
2\int_0^\sigma (X_t-X^n_t)\mathrm{d}M_t.
\end{equation*}
Let us show $X^n \in  \mathcal{T}_u(S,\beta,\sigma)$.
The convergence  (\ref{sup}) is apparent by definition.
Since
\begin{equation} \label{costrep}
\begin{split}
C[S,\beta;X^n]_{\sigma} =& 
\sum_{0 < t \leq \sigma} S_tK_t|\Delta X^n_t|^{\beta - 2}|\Delta
X^n_t|^2 1_{\{|\Delta X^n_t| > 0\}}\\
=& \epsilon_n^{\beta -2} \sum_{ j \geq 1, \tau^n_j \leq \sigma} 
S_{\tau^n_j}K_{\tau^n_j}{\hat{S}_{\tau^n_{j-1}}^{\beta-2}} 
|X_{\tau^n_j} - X_{\tau^n_{j-1}}|^2,
\end{split}
\end{equation}
there exists a constant $c >0 $ such that
\begin{equation*}
\frac{1}{c} 
\mathbb{E}\left[
\sum_{j=0}^{\infty}
|X_{\tau^n_{j+1} \wedge \sigma} - X_{\tau^n_{j} \wedge \sigma}|^2
\right] \leq \epsilon_n^{2-\beta} 
\mathbb{E}[C[S,\beta;X^n]_{\sigma}] \leq  
c 
\mathbb{E}\left[
\sum_{j=0}^{\infty}
|X_{\tau^n_{j+1} \wedge \sigma} - X_{\tau^n_{j} \wedge \sigma}|^2
\right].
\end{equation*}
Since
\begin{equation*}
\mathbb{E}\left[
\sum_{j=0}^{\infty}
|X_{\tau^n_{j+1} \wedge \sigma} - X_{\tau^n_{j} \wedge \sigma}|^2
\right] = 
\mathbb{E}[\langle X \rangle_\sigma]
+ 
 \mathbb{E}[\int_0^{\sigma}(X_s-X^n_s)H_s\mathrm{d}\langle X
\rangle_s],
\end{equation*}
we obtain $\mathbb{E}[C[S,\beta;X^n]_{\sigma}]= O(\epsilon_n^{\beta - 2})$.
On the other hand, 
\begin{equation*}
\langle Z[X^n]\rangle_{\sigma} \leq \sup_{ t \in [0,\sigma]}
\{K_t|X_t-X^n_t|^2\} \langle X \rangle_{\sigma}
\leq \epsilon_n^2 \langle X \rangle_{\sigma}
\sup_{ t \in [0,\sigma]}\hat{S}_t
\sup_{ t \in [0,\sigma]}K_t,
\end{equation*}
and so, we conclude that
\begin{equation}\label{cz}
|\mathbb{E}[C[S,\beta;X^n]_{\sigma}]|^{2/(2-\beta)}
\langle Z[X^n]\rangle_{\sigma}
\end{equation}
is uniformly integrable. 
So far we showed that 
$X^n \in \mathcal{T}_u(S,\beta,\sigma)$.
The uniform integrability of
\begin{equation*}
\frac{C[S,\beta;X^n]_{\sigma}}
{\mathbb{E}[C[S,\beta;X^n]_{\sigma}]}
\end{equation*}
also follows from that of (\ref{qv}) in light of (\ref{costrep}).
With the aid of the uniform integrability of (\ref{qv}) and (\ref{cz}),
repeating the same argument as in the proof of Theorem~\ref{R1},
the convergence result follows from the fact that
\begin{equation*}
\begin{split}
\epsilon_n^{2-\beta}
C[S,\beta;X^n]_{\sigma} =& 
\sum_{ j \geq 1, \tau^n_j \leq \sigma} 
S_{\tau^n_j}K_{\tau^n_j} \hat{S}_{\tau^n_{j-1}}^{\beta-2}
|X_{\tau^n_j} - X_{\tau^n_{j-1}}|^2
\\&\to \int_0^\sigma S_t \hat{S}_t^{\beta-2} K_t \mathrm{d}\langle X \rangle_t
 = (S \hat{S}^{\beta-2}) \cdot \langle Y \rangle_\sigma,
\\
\epsilon_n^{-2} 
\sum_{ j \geq 1, \tau^n_j \leq \sigma} 
K_{\tau^n_{j-1}}(X_{\tau^n_{j}} - 
X_{\tau^n_{j-1}})^4 = & 
\sum_{ j \geq 1, \tau^n_j \leq \sigma} 
K_{\tau^n_{j-1}} \hat{S}_{\tau^n_{j-1}}^{2} 
|X_{\tau^n_j} - X_{\tau^n_{j-1}}|^2
\\ &
\to \int_0^\sigma \hat{S}^{2}_t K_t \mathrm{d}\langle X \rangle_t
= \hat{S}^{2} \cdot \langle Y \rangle_\sigma
\end{split}
\end{equation*}
in probability as $n\to \infty$. 
\hfill////
\begin{rem} \upshape
The assumption $\langle Y \rangle = K \cdot \langle X \rangle$ implies
 in particular that $Y$ is quasi-left-continuous.
That $Y$ is quasi-left-continuous is equivalent to that
$Y$ has no predictable jump time.
See Jacod and Shiryaev~\cite{JS} for more details.
For example, the L\'evy processes are quasi-left-continuous.
Of course so are the continuous semimartingales.
The asymptotic efficiency of (\ref{eff}) is no more true if $Y$ is not
 quasi-left continuous. In fact, if there is a predictable time $\tau$
such that $Y_{\tau} \neq Y_{\tau-}$, it is apparently more efficient to include
 $\tau$, or more precisely, a time immediately before $\tau$ 
into the sequence of stopping times for discretization. 
This is possible because $\tau$ is  predictable.
\end{rem}

\section{Efficiency for possibly biased Riemann sums}
\subsection{The case of $\beta \in [0,1]$}
The class $\mathcal{T}_u(S,\beta,\sigma)$
was a set of unbiased schemes, that is, 
$\{X^n\}$ of the form $X^n = X[\tau^n]$, $\tau^n \in
\mathcal{T}$.
As an approximating sequence $X^n$ 
to $X$, we may consider more general 
simple predictable processes.
In this section we answer the question that the scheme (\ref{eff})
is asymptotically efficient in a more general class of 
simple predictable processes or not. 
First we get a positive answer for $\beta \in [0,1]$.
The result improves Fukasawa~\cite{SAFA} for the case
$\beta = 0$.
Denote by $\mathcal{T}(S,0,\sigma)$ the set of the  sequences 
$X^n$ of simple predictable processes such that
that there exists a sequence of stopping times $\sigma^m$ with 
$\sigma^m \to \infty$ as $m \to \infty$, 
\begin{enumerate}
\item for each $m$,
\begin{equation*}
\sup_{t \in [0,\sigma^m]}|X^n_t-X_t|
\end{equation*}
is uniformly bounded and converges to $0$
in probability as $n \to \infty$,  and
\item for each $m$,
\begin{equation*}
\mathbb{E}[C[S,0;X^n]_{\sigma^m}]
\langle Z[X^n]\rangle_{\sigma^m}
\end{equation*}
is uniformly integrable in $n$.
\end{enumerate}
For $\beta \in (0,2)$, we need additional conditions from technical
point of view.
We define $\mathcal{T}(S,\beta,\sigma)$ 
for $\beta \in (0,2)$ as
the set of the sequences 
$X^n$ of simple predictable processes such 
that there exists a sequence of stopping times $\sigma^m$ with 
$\sigma^m \to \infty$ as $m \to \infty$, 
\begin{enumerate}
\item for each $m$,
\begin{equation*}
\sup_{t \in [0,\sigma^m]}|X^n_t-X_t|,\ \ 
\sup_{t \in [0,\sigma^m]}\left|
\frac{\Delta X^n_t}{\Delta X[\tau[X^n]]_t}-1 \right| 
\end{equation*}
are uniformly bounded and converge to $0$
in probability as $n \to \infty$, where $0/0$ is understood as $1$, and
\item for each $m$,
\begin{equation*}
|\mathbb{E}[C[S,\beta;X^n]_{\sigma^m}]|^{2/(2-\beta)}
\langle Z[X^n]\rangle_{\sigma^m}, \ \ 
\frac{C[S,\beta;X[\tau[X^n]]]_{\sigma^m}}
{\mathbb{E}[C[S,\beta;X[\tau[X^n]]]_{\sigma^m}]}
\end{equation*}
are uniformly integrable in $n$.
\end{enumerate}
The convergence of the ratio between $\Delta X[\tau[X^n]]$ and
$\Delta X^n$ to $1$ means that $X^n$ cannot be too biased.
Of course it always holds if $X^n$ is unbiased since 
$X[\tau[X^n]] = X^n$.
The uniform integrability of the normalized 
cost function associated with $X[\tau[X^n]]$ is  reasonable in that
it requires the sequence of stopping times $\tau[X^n]$
to be sufficiently regular.
By Theorem~\ref{R2},
the scheme $\{X^n\}$ defined by (\ref{hitting})
is an element of  $\mathcal{T}(S,\beta,\sigma)$ for 
any $\beta \in [0,2)$.
Therefore, the following theorem asserts that
the scheme $\{X^n\}$ defined by (\ref{eff})
is asymptotically efficient in the class $\mathcal{T}(S,\beta,\sigma)$
 if $\beta \in [0,1]$.
\begin{thm}  \label{main1}
Let $\beta \in [0,1]$.
The inequality (\ref{cx}) holds for 
 all $\{X^n\} \in \mathcal{T}(S,\beta,\sigma)$.
\end{thm}
{\it Proof: } Write $\tau^n = \tau[X^n]$ for brevity.
By the usual localization procedure, 
we may and do suppose without loss of generality that
$X,  H \cdot M,  \langle X \rangle, K, 1/K, S, 1/S$ and  $H$
are bounded up to $\sigma$, that $\sup_{t \in [0,\sigma]}|X^n_t-X_t|$
is uniformly bounded and converge to $0$,  and that
$|\mathbb{E}[C[S,\beta;X[\tau[X^n]]]_{\sigma}]|^{2/(2-\beta)}
\langle Z[X^n]\rangle_{\sigma}$
is uniformly integrable in $n$. 
For the case $\beta \in (0,1]$, we may have additionally
that
\begin{equation*}
\sup_{t \in [0,\sigma]}\left|
\frac{\Delta X^n_t}{\Delta X[\tau^n]_t}-1 \right| 
\end{equation*}
is uniformly bounded and converge to $0$,  and that
\begin{equation*}
\frac{C[S,\beta;X[\tau^n]]_{\sigma}}
{\mathbb{E}[C[S,\beta;X[\tau^n]]_{\sigma}]}
\end{equation*}
is  uniformly integrable in $n$. 
Define $K$ and  $K[\tau^n]$ 
as in the proof of Theorem~\ref{R1}.
By It$\hat{\text{o}}$'s formula,
\begin{equation*}
\begin{split}
\langle Z[X^n] \rangle_\sigma \geq
& \int_0^\sigma |X_s - X^n_s|^2
 K[\tau^n]_s\mathrm{d}\langle X \rangle_s
+ \int_0^\sigma |X_s - X^n_s|^2(K_s - K[\tau^n]_s)\mathrm{d}\langle X
 \rangle_s 
\\ =& \frac{1}{6} \sum_{j=0}^{\infty}K_{\tau^n_j}
((\Delta_j  + \delta_j)^4 - \delta_j^4 ) 
- \frac{2}{3}\int_0^\sigma K[\tau^n]_s(X_s-X^n_s)^3\mathrm{d}X_s \\
& + 
\int_0^\sigma |X_s - X^n_s|^2(K_s - K[\tau^n]_s)\mathrm{d}\langle X
 \rangle_s,
\end{split}
\end{equation*}
where 
$\Delta_j = X_{\tau^n_{j+1} \wedge \sigma} - X_{\tau^n_j \wedge \sigma}$
and  
$\delta_j = X_{\tau^n_j \wedge \sigma} - X^n_{\tau^n_j \wedge \sigma}$.
As before, we can show that
\begin{equation*}
\begin{split}
&\lim_{n\to \infty}
|\mathbb{E}[C[S,\beta;X^n]_{\sigma}]|^{2/(2-\beta)}
\mathbb{E}[
\int_0^\sigma |X_s - X^n_s|^2(K_s - K[\tau^n]_s)\mathrm{d}\langle X
 \rangle_s] = 0, \\
& \lim_{n\to \infty}
|\mathbb{E}[C[S,\beta;X^n]_{\sigma}]|^{2/(2-\beta)}
\mathbb{E}[\int_0^\sigma K[\tau^n]_s(X_s-X^n_s)^3\mathrm{d}X_s
] = 0
\end{split}
\end{equation*}
by the uniform integrability with the aid of Lemma~\ref{supsup}. Put
\begin{equation*}
F_t = \exp \left\{
\int_0^t H_s \mathrm{d}M_s - \frac{1}{2}\int_0^t H_s^2 \mathrm{d}\langle
M \rangle_s
\right\}.
\end{equation*}
Since
\begin{equation*}
\mathbb{E}[F_{\tau^n_{j+1}}/F_{\tau^n_j}] = 1, \ \ 
\sup_{t \geq 0, j \geq 0} \left|
1 - \frac{F_{t \wedge \tau^n_{j+1} \wedge \sigma}}
{F_{\tau^n_j \wedge \sigma}}
\right| \to 0
\end{equation*}
in probability, 
again by It$\hat{\text{o}}$'s formula, we have that
\begin{equation*}
\begin{split}
&\liminf_{n\to \infty}
|\mathbb{E}[C[S,\beta;X^n]_{\sigma}]|^{2/(2-\beta)}
\mathbb{E}[
\langle Z[X^n] \rangle_\sigma]\\
& = \frac{1}{6}\liminf_{n\to \infty}
|\mathbb{E}[C[S,\beta;X^n]_{\sigma}]|^{2/(2-\beta)}
\mathbb{E}\left[
\sum_{j=0}^{\infty}K_{\tau^n_j} 
((\Delta_j  + \delta_j)^4 - \delta_j^4 )\right]
\\
& = \frac{1}{6}\liminf_{n\to \infty}
|\mathbb{E}[C[S,\beta;X^n]_{\sigma}]|^{2/(2-\beta)}
\mathbb{E}\left[
\sum_{j=0}^{\infty}K_{\tau^n_j} 
((\Delta_j  + \delta_j)^4 - \delta_j^4 )
\frac{F_{\tau^n_{j+1} \wedge \sigma}}{F_{\tau^n_j \wedge
 \sigma}}
\right] 
\\
& = \frac{1}{6} \liminf_{n\to \infty}
|\mathbb{E}[C[S,\beta;X^n]_{\sigma}]|^{2/(2-\beta)}
\mathbb{E}\left[
\sum_{j=0}^{\infty}K_{\tau^n_j} \hat{\mathbb{E}}_j[
((\Delta_j  + \delta_j)^4 - \delta_j^4 )]
\right], 
\end{split}
\end{equation*}
where $\hat{\mathbb{E}}_j[A]$ refers to
the conditional expectation 
$\mathbb{E}[AF_{\tau^n_{j+1} \wedge \sigma}/F_{\tau^n_j \wedge
 \sigma}| \mathcal{F}_{\tau^n_j \wedge \sigma}]$ 
for a random variable $A$.
Notice that under $\hat{\mathbb{E}}_j$, 
$X_{t \wedge \tau^n_{j+1}\wedge \sigma} -
X_{t \wedge \tau^n_j \wedge \sigma}$ 
is a martingale. Therefore, 
\begin{equation*}
\begin{split}
\hat{\mathbb{E}}_j[
((\Delta_j  + \delta_j)^4 - \delta_j^4 )]
=& \hat{\mathbb{E}}_j[\Delta_j^4]
+4 \delta_j\hat{\mathbb{E}}_j[\Delta_j^3]
+ 6\delta_j^2\hat{\mathbb{E}}_j[\Delta_j^2]
\\
=& 
6\hat{\mathbb{E}}_j[\Delta_j^2] \left(
\delta_j +
 \frac{1}{3} \frac{\hat{\mathbb{E}}_j[\Delta_j^3]}
{\hat{\mathbb{E}}_j[\Delta_j^2]}\right)^2
+ \hat{\mathbb{E}}_j[\Delta_j^4] - \frac{2}{3}
 \frac{|\hat{\mathbb{E}}_j[\Delta_j^3]|^2}
{\hat{\mathbb{E}}_j[\Delta_j^2]} \\
\geq & 
 \frac{|\hat{\mathbb{E}}_j[\Delta_j^2]|^{(4-\beta)/(2-\beta)}}
{|\hat{\mathbb{E}}_j[|\Delta_j|^{\beta}]|^{2/(2-\beta)}}.
\end{split}
\end{equation*}
Here, we have used Lemma~\ref{KS1} for $\beta \in [0,1) $
and (\ref{Fukasawa}) for $\beta = 1$.

By H$\ddot{\text{o}}$lder's inequality,
\begin{equation*}
\begin{split}
&\mathbb{E}
\left[
\sum_{\tau^n_{j+1}\leq \sigma}^{\infty}
|S_{\tau^n_j}|^{2/(4-\beta)}K_{\tau^n_j} 
\hat{\mathbb{E}}_j[\Delta_j^2]
\right] \\ 
&
\leq 
\left|
\mathbb{E}
\left[
\sum_{j=0}^{\infty} 
\frac{K_{\tau^n_j}|\hat{\mathbb{E}}_j[\Delta_j^2]|^p}
{|\hat{\mathbb{E}}_j[|\Delta_j|^{\beta}]|^{2/(2-\beta)}}
\right]
\right|^{1/p}
\left|
\mathbb{E}
\left[
\sum_{\tau^n_{j+1} \leq \sigma} S_{\tau^n_j} K_{\tau^n_j}
\hat{\mathbb{E}}_j[|\Delta_j|^{\beta}]
\right]
\right|^{1/q},
\end{split}
\end{equation*}
where $p = (4-\beta)/(2-\beta)$ and 
$q = p/(p-1) = (4-\beta)/2$.
Since
\begin{equation*}
\sup_{t \in [0,\sigma]}\left|
\frac{|\Delta X^n_t|^\beta}{|\Delta X[\tau^n]_t|^\beta}-1 \right|, \ \ 
\sup_{\tau^n_{j+1} \leq \sigma}\left|
\frac{S_{\tau^n_j}K_{\tau^n_j}}{S_{\tau^n_{j+1}}K_{\tau^n_{j+1}}}
\frac{F_{\tau^n_{j+1} }}
{F_{\tau^n_{j}}} -1
\right|
\end{equation*}
are uniformly bounded and converge to $0$ in probability, we get
\begin{equation*}
\lim_{n \to \infty}
\frac{\mathbb{E}
\left[
\sum_{\tau^n_{j+1} \leq \sigma} S_{\tau^n_j} K_{\tau^n_j}
\hat{\mathbb{E}}_j[|\Delta_j|^{\beta}]
\right]}{\mathbb{E}[C[S,\beta;X[\tau^n]]_\sigma]} =
 \lim_{n \to \infty}
\frac{\mathbb{E}[C[S,\beta;X[\tau^n]]_\sigma]}
{\mathbb{E}[C[S,\beta;X^n]_\sigma]} = 1.
\end{equation*}
Here we have used the 
uniform integrability of
$C[S,\beta;X[\tau^n]]_{\sigma}/ \mathbb{E}[C[S,\beta;X[\tau^n]]_{\sigma}]$ 
for $\beta \in (0,1]$.
This is trivial if $\beta = 0$.

By the bounded convergence theorem, 
\begin{equation*}
\lim_{n\to \infty}
\mathbb{E}
\left[
\sum_{\tau^n_{j+1}\leq \sigma}
|S_{\tau^n_j}|^{2/(4-\beta)}K_{\tau^n_j} 
\hat{\mathbb{E}}_j[\Delta_j^2]
\right]
= \mathbb{E}[(S^{2/(4-\beta)} K)\cdot \langle X\rangle_{\sigma}],
\end{equation*}
which completes the proof.
\hfill////

\subsection{The case of $\beta \in (1,2)$}
Here we show that the unbiased scheme $X^n$ 
defined by (\ref{eff}) is no more efficient for 
$\beta \in (1,2)$.
We give a lower bound which is one third the previous one and 
construct a biased scheme which asymptotically attains it.
\begin{thm}
Let $\beta \in (1,2)$. 
For all 
$\{X^n\} \in \mathcal{T}(S,\beta,\sigma)$,
\begin{equation*}
\liminf_{n \to \infty}
|\mathbb{E}[C[S,\beta;X^n]_\sigma]|^{2/(2-\beta)}
\mathbb{E}[\langle Z[X^n]\rangle_\sigma]
\geq \frac{1}{18}
|\mathbb{E}[
(S^{2/(4-\beta)}K)  \cdot 
\langle X \rangle_\sigma]|^{(4-\beta)/(2-\beta)}.
\end{equation*}
\end{thm}
{\it Proof: } Just use Lemma~\ref{KS2} with $\alpha = 2/3$
instead of Lemma~\ref{KS1} in the proof of Theorem~\ref{main1}.
The rest is the same.
\hfill////
\begin{thm}
Suppose that $\langle Y \rangle = K \cdot \langle X \rangle$.
Let $\beta \in (1,2)$ and
 $\epsilon_n$ be a positive sequence with $\epsilon_n \to 0$ as 
$n\to  \infty$. For $\gamma  \in \mathbb{R}$, 
define $\tau^n(\gamma) = \{\tau^n_j(\gamma)\}$ as
\begin{equation}\label{gamma}
\begin{split}
&\tau^n_0(\gamma) = 0, \ \ \tau^n_{j+1}(\gamma) =
\min\{\tau^n_{j+1}(\gamma, +), \tau^n_{j+1}(\gamma, -) \}, \\ 
& \tau^n_{j+1}(\gamma, +) = 
\inf\left\{t > \tau^n_j(\gamma); X_t - X_{\tau^n_j(\gamma)} \geq 
\epsilon_n e^\gamma S_{\tau^n_j(\gamma)}^{1/(4-\beta)} \right\}, \\
& \tau^n_{j+1}(\gamma, -) = 
\inf\left\{t > \tau^n_j(\gamma);
X_t - X_{\tau^n_j(\gamma)} \leq 
\epsilon_n e^{-\gamma} S_{\tau^n_j(\gamma)}^{1/(4-\beta)} \right\}.
\end{split}
\end{equation}
Define a sequence of simple predictable processes $X^n(\gamma)$ as
\begin{equation} \label{Xgamma}
X^n(\gamma) = X[\tau^n(\gamma)] + \frac{2}{3} \epsilon_n
\sinh(\gamma)  S[\tau^n(\gamma)]^{1/(4-\beta)},
\end{equation}
where $S[\tau^n(\gamma)]_t = S_{\tau^n_j(\gamma)}$ for 
$t \in  [\tau^n_j(\gamma), \tau^n_{j+1}(\gamma))$.
Then $\{X^n(\gamma)\} \in \mathcal{T}(S,\beta,\sigma)$. 
Moreover if
$X,\langle X \rangle, H\cdot M, H,  K,   1/K,  S$ and  $1/S$ 
are bounded up to $\sigma$,
then
\begin{equation*}
\lim_{n \to \infty} 
|\mathbb{E}[C[S,\beta;X^n(\gamma)]_{\sigma}]|^{2/(2-\beta)}
\mathbb{E}[\langle Z[X^n(\gamma)]\rangle_{\sigma}]
= \frac{F(|\gamma|)}{6}
|\mathbb{E}[
S^{2/(4-\beta)}  \cdot 
\langle Y \rangle_{\sigma}]|^{(4-\beta)/(2-\beta)},
\end{equation*}
where $F$ is a continuous function
with $F(0) = 1$ and $F(\infty) = 1/3$. More explicitly,
\begin{equation*}
F(x) = F(x,\beta) = \frac{ 4|\cosh(x)|^2 -1}
{3|\cosh(x)|^{2/(2-\beta)} |\cosh((\beta-1)x)|^{-2/(2-\beta)} }.
\end{equation*}
\end{thm}
{\it Proof: }
By the usual localization procedure, 
we may and do suppose without loss of generality that
$X,  H \cdot M,  \langle X \rangle, K, 1/K, S, 1/S$ and  $H$
are bounded up to $\sigma$.
Put $X^n = X^n(\gamma)$ and $\tau^n = \tau[X^n] = \tau^n(\gamma)$
for brevity.
Then it follows from definition that
\begin{equation*}
\sup_{t \in [0,\sigma]}|X^n_t-X_t|, \ \ 
\sup_{t \in [0,\sigma]}\left|
\frac{\Delta X^n_t}{\Delta X[\tau^n]_t}-1 \right| 
\end{equation*}
are uniformly bounded and converge to $0$.
By the same argument as in the proof of Theorem~\ref{R2},
we have that
\begin{equation*}
|\mathbb{E}[C[S,\beta;X[\tau^n]]_{\sigma}]|^{2/(2-\beta)}
\langle Z[X^n]\rangle_{\sigma}, \ \ 
\frac{C[S,\beta;X[\tau^n]]_{\sigma}}
{\mathbb{E}[C[S,\beta;X[\tau^n]]_{\sigma}]}
\end{equation*}
are uniformly integrable in $n$.
Since these imply in particular that
\begin{equation*}
\lim_{n\to \infty}\frac{\mathbb{E}[C[S,\beta;X^n]_{\sigma}]}
{\mathbb{E}[C[S,\beta;X[\tau^n]]_{\sigma}]}  = 1,
\end{equation*}
we conclude  $\{X^n\} \in \mathcal{T}(S,\beta,\sigma)$.

Let $\Delta_j = X_{\tau^n_{j+1}} - X_{\tau^n_j}$ and
$\delta_j = X_{\tau^n_j} - X^n_{\tau^n_j}$.
Then we obtain, in a similar manner to the proof of Theorem~\ref{main1}, that
\begin{equation*}
\begin{split}
&\lim_{n\to \infty}
|\mathbb{E}[C[S,\beta;X^n]_{\sigma}]|^{2/(2-\beta)}
\mathbb{E}[
\langle Z[X^n] \rangle_\sigma]\\
& = \frac{1}{6} \liminf_{n\to \infty}
|\mathbb{E}[C[S,\beta;X^n]_{\sigma}]|^{2/(2-\beta)}
\mathbb{E}\left[
\sum_{\tau^n_{j} \leq \sigma} K_{\tau^n_j} \hat{\mathbb{E}}_j[
((\Delta_j  + \delta_j)^4 - \delta_j^4 )]
\right]
\end{split}
\end{equation*}
and that
\begin{equation*}
\hat{\mathbb{E}}_j[
((\Delta_j  + \delta_j)^4 - \delta_j^4 )]
= 
6\hat{\mathbb{E}}_j[\Delta_j^2] \left(
 \delta_j +
 \frac{1}{3} \frac{\hat{\mathbb{E}}_j[\Delta_j^3]}
{\hat{\mathbb{E}}_j[\Delta_j^2]}\right)^2
+ \hat{\mathbb{E}}_j[\Delta_j^4] - \frac{2}{3}
 \frac{|\hat{\mathbb{E}}_j[\Delta_j^3]|^2}
{\hat{\mathbb{E}}_j[\Delta_j^2]},
\end{equation*}
where $\hat{\mathbb{E}}_j[A]$ refers to
the conditional expectation 
$\mathbb{E}[AF_{\tau^n_{j+1}}/F_{\tau^n_j }| \mathcal{F}_{\tau^n_j}]$ 
for a random variable $A$.
By the optional sampling theorem,
\begin{equation*}
\hat{\mathbb{E}}_j[I\{\Delta_j = \epsilon_n e^\gamma
 S_{\tau^n_j}^{1/(4-\beta)}\}]
 = \frac{e^{-\gamma}}{ e^{\gamma} + e^{-\gamma}}, \ \ 
\hat{\mathbb{E}}_j[I\{\Delta_j = -\epsilon_n e^{-\gamma}
S_{\tau^n_j}^{1/(4-\beta)}\}]
 = \frac{e^{\gamma}}{ e^{\gamma} + e^{-\gamma}}
\end{equation*}
and so, 
\begin{equation*}
\begin{split}
&\hat{\mathbb{E}}_j[\Delta_j] = 0, \ \ 
\hat{\mathbb{E}}_j[\Delta_j^2] = 
\epsilon_n^2  S_{\tau^n_j}^{2/(4-\beta)}, \ \
\hat{\mathbb{E}}_j[\Delta_j^3] = 
2\epsilon_n^3 \sinh(\gamma) S_{\tau^n_j}^{3/(4-\beta)},\\
&\hat{\mathbb{E}}_j[|\Delta_j|^\beta] = 
\epsilon_n^\beta  \frac{\cosh((\beta-1)\gamma)}{\cosh(\gamma)} 
S_{\tau^n_j}^{\beta/(4-\beta)}.
\end{split}
\end{equation*} 
Moreover by Lemma~\ref{KS2},
\begin{equation*}
\begin{split}
&\hat{\mathbb{E}}_j[\Delta_j^4]- \frac{2}{3}
\frac{|\hat{\mathbb{E}}_j[\Delta_j^3]|^2}{\hat{\mathbb{E}}_j[\Delta_j^2]}  \\
&= F(|\gamma|)
\frac{|\hat{\mathbb{E}}_j[\Delta_j^2]|^{(4-\beta)/(2-\beta)}}
{|\hat{\mathbb{E}}_j[|\Delta_j|^{\beta}]|^{2/(2-\beta)}}
\\&= F(|\gamma|) \left|\frac{\cosh((\beta-1)\gamma)}{\cosh(\gamma)} \right|^{-2/(2-\beta)}
\epsilon_n^2 S_{\tau^n_j}^{2/(4-\beta)} \hat{\mathbb{E}}_j[\Delta_j^2]
\end{split}
\end{equation*}
with $F = F(\cdot,\beta)$, which satisfies
$F(|\gamma|) \to 1/3$ as $|\gamma| \to \infty$.
By definition of $X^n$, we have
\begin{equation*}
 \delta_j +
 \frac{1}{3} \frac{\hat{\mathbb{E}}_j[\Delta_j^3]}
{\hat{\mathbb{E}}_j[\Delta_j^2]} = 0.
\end{equation*}
Therefore,
\begin{equation*}
\lim_{n\to \infty}
\epsilon_n^{-2}
\mathbb{E}[
\langle Z[X^n] \rangle_\sigma] = \frac{1}{6} F(|\gamma|)
 \left|\frac{\cosh((\beta-1)\gamma)}{\cosh(\gamma)}  \right|^{-2/(2-\beta)}
\mathbb{E}[S^{2/(4-\beta)} \cdot \langle Y \rangle_{\sigma}].
\end{equation*}
On the other hand,
\begin{equation*}
\begin{split}
&\lim_{n\to \infty}
\epsilon_n^{2-\beta}\mathbb{E}[C[S,\beta;X[\tau^n]_{\sigma}] \\
&= \lim_{n\to \infty} \mathbb{E}\left[
\sum_{\tau^n_{j}\leq \sigma}
S_{\tau^n_{j}} K_{\tau^n_j} \hat{\mathbb{E}}_j[|\Delta_j|^{\beta}]
\right]\\
&= \frac{\cosh((\beta-1)\gamma)}{\cosh(\gamma)}
\lim_{n\to \infty} \mathbb{E}\left[
\sum_{\tau^n_{j}\leq \sigma}
S_{\tau^n_{j}}^{2/(4-\beta)} K_{\tau^n_j} \hat{\mathbb{E}}_j[\Delta_j^2]
\right] \\
&= \frac{\cosh((\beta-1)\gamma)}{\cosh(\gamma)}
 \mathbb{E}[S^{2/(4-\beta)} \cdot \langle Y \rangle_{\sigma}].
\end{split}
\end{equation*}
These convergences give the result. \hfill////

\begin{rem} \upshape
The use of hitting times is essential to have a good performance.
In fact if we consider a class of simple predictable processes $X^n$
such that $\tau[X^n]_{j+1} - \tau[X^n]_j$ is
 $\mathcal{F}_{\tau[X^n]_j}$-measurable for each $j \geq 0$,
then we can show that
\begin{equation*}
\liminf_{n \to \infty}
|\mathbb{E}[C[S,\beta;X^n]_\sigma]|^{2/(2-\beta)}
\mathbb{E}[\langle Z[X^n]\rangle_\sigma]
\geq \frac{1}{2}
|\mathbb{E}[
(S^{2/(4-\beta)}K)  \cdot 
\langle X \rangle_\sigma]|^{(4-\beta)/(2-\beta)}
\end{equation*}
when, for example, $X = Y$ and it is a Brownian motion.
This is because the kurtosis 
$\hat{\mathbb{E}}_j[\Delta_j^4] | \hat{\mathbb{E}}_j[\Delta_j^2]|^{-2}$
and skewness
$\hat{\mathbb{E}}_j[\Delta_j^3] | \hat{\mathbb{E}}_j[\Delta_j^2]|^{-3/2}$
of a conditionally standard normal random variable $\Delta_j$ are
 $3$ and $0$ respectively, while the lower bound of kurtosis
is $1$ attained by Bernoulli random variables.
The above measurability condition was supposed 
in Genon-Catalot and Jacod~\cite{GJ}.
\end{rem}

\section{Exponential utility maximization}
The schemes $X^n = X[\tau^n]$ with (\ref{eff}) and $X^n = X^n(\gamma)$ defined by
(\ref{Xgamma}) with (\ref{gamma}) are efficient for 
$\beta \in [0,1]$ and $\beta \in (1,2)$ respectively in that they
attain the asymptotic lower bound of
\begin{equation*}
|\mathbb{E}[C[S,\beta;X^n]_{\sigma}]|^{2/(2-\beta)}
\mathbb{E}[\langle Z[X^n]\rangle_{\sigma}]
\end{equation*}
for a reasonable class of approximating simple predictable processes $X^n$.
In the financial context of discrete hedging, we may interpret 
the cost function $C[S,\beta;\hat{X}]_{\sigma}$ 
as the cumulative transaction cost associated to the rebalancing scheme
$\hat{X}$. If we do so, then a more natural criterion for the optimality of
$\hat{X}$ should be given in terms of the expected utility of the terminal
wealth $-Z[\hat{X}]_{\sigma}-  C[S,\beta;\hat{X}]_{\sigma}$. 
In this section, we see that the efficient schemes maximize a scaling
limit of the exponential utility
\begin{equation*}
1- \mathbb{E}[\exp\{- \alpha_n(-Z[X^n]_{\sigma}- 
C[S^n,\beta;X^n]_{\sigma})\}], \ \ 
S^n = \kappa_n S, \ \ \alpha_n \to \infty,\ \ 
\alpha_n \kappa_n \to 0.
\end{equation*}
Here $\kappa_n$ is a deterministic sequence, which we interpret as the
coefficient of the transaction costs.
Letting $\kappa_n \to 0$, we try to obtain an asymptotic but explicit
solution for the maximization problem which can be expected to
 have a good performance when $\kappa_n$  is sufficiently small.
If $\kappa_n \to 0$, then we can make both 
$\langle Z[X^n] \rangle_\sigma$ and $C[S^n,\beta;X^n]_{\sigma}$ converge
to $0$ by taking any $\{X^n\} \in \mathcal{T}(S,\beta,\sigma)$ 
such that 
 $\sup_{t \in [0,\sigma]}|X^n_t-X_t| \to 0$ sufficiently slow. 
To find effective $X^n$ among others,  we 
consider a scaling limit by letting $\alpha_n$, the
risk-aversion parameter, diverge.
In this section we assume $Y$ to be continuous in addition.
By Jacod's theorem of stable convergence of semimartingales,
if there exists a continuous process $V$ such that
\begin{equation} \label{jacod}
\alpha_n^2 \langle Z[X^n] \rangle_t \to V_t, \ \ 
\alpha_n \langle Z[X^n],Y \rangle_t \to 0
\end{equation}
in probability for all $t \geq 0$, then $\alpha_n Z[X^n]$ converges
$\mathcal{F}$-stably in law to a time-changed Brownian motion $W_V$,
where $W$ is a standard Brownian motion which is independent of
$\mathcal{F}$. See Fukasawa~\cite{AAP} for more details and 
sufficient conditions for (\ref{jacod}). 
Note that the second condition of (\ref{jacod}) is to make the
replication error $Z[X^n]$ asymptotically neutral to the market return. 
If in addition $\alpha_n C[S^n,\beta;X^n]_{\sigma}$ converges to 
a random variable $C_\sigma$ in probability, then
\begin{equation*}
\alpha_nZ[X^n]_{\sigma} + \alpha_n 
C[S^n,\beta;X^n]_{\sigma} \to W_{V_{\sigma}} + C_\sigma
\end{equation*}
in law. The limit law is a mixed normal distribution with
conditional mean $C_\sigma$ and conditional variance $V_{\sigma}$.
This implies in particular that 
\begin{equation*}
1- \mathbb{E}[\exp\{- \alpha_n(-Z[X^n]_{\sigma}- 
C[S^n,\beta;X^n]_{\sigma})\}] \to
1- \mathbb{E}[\exp\{C_\sigma + \frac{1}{2}V_{\sigma}\}]
\end{equation*}
under the uniform integrability condition on
$\exp\{\alpha_n(Z[X^n]_{\sigma}+ C[S^n,\beta;X^n]_{\sigma})\}$.
Then the maximization of the exponential utility reduces to the minimization
of $C_\sigma + V_{\sigma}/2$. Under the additional assumptions that
\begin{equation*}
\alpha_n^4 \sum_{j=0}^{\infty}
\mathbb{E}[|\langle X \rangle_{\tau^n_{j+1} \wedge \sigma} - \langle X
\rangle_{\tau^n_{j} \wedge \sigma}|^4
| \mathcal{F}_{\tau^n_{j} \wedge \sigma}] \to 0
\end{equation*}
in probability with $\tau^n = \tau[X^n]$ and that
\begin{equation*}
\alpha_n^{(6-2\beta)/(2-\beta)} \kappa_n^{2/(2-\beta)} \to \mu > 0,
\end{equation*}
we obtain that
\begin{equation*}
C_\sigma^{2/(2-\beta)}V_{\sigma} \geq \frac{\mu}{6} | S^{2/(4-\beta)}
 \cdot \langle Y \rangle_\sigma |^{(4-\beta)/(2-\beta)}
\end{equation*}
for $\beta \in [0,1]$ by a similar argument to the proof of
Theorem~\ref{main1} with the aid of Lemma~A.2 of Fukasawa~\cite{AAP}.
This is in fact an extension of Theorems~2.7 and 2.8 of
Fukasawa~\cite{AAP}.
It follows then that
\begin{equation*}
\begin{split}
C_\sigma + \frac{1}{2}V_{\sigma} &\geq 
C_\sigma + \frac{\mu}{12} | S^{2/(4-\beta)}
 \cdot \langle Y \rangle_\sigma |^{(4-\beta)/(2-\beta)}
C_\sigma^{-2/(2-\beta)}  \\ &\geq
\hat{\mu}
S^{2/(4-\beta)}
 \cdot \langle Y \rangle_\sigma,
\end{split}
\end{equation*}
where
\begin{equation*}
\hat{\mu} = \left|\frac{\mu}{6(2-\beta)} \right|^{(2-\beta)/(4-\beta)}  +
 \frac{\mu}{12}
\left|\frac{\mu}{6(2-\beta)} \right|^{-2/(4-\beta)}.
\end{equation*}
Here we have used the fact that for given $c > 0$, 
$\min_{x > 0}\{x + c x^{-2/(2-\beta)}\}$ is attained at
$x = (2c/(2-\beta))^{(2-\beta)/(4-\beta)}$.
Therefore,
\begin{equation*}
\lim_{n \to \infty}
\{1- \mathbb{E}[\exp\{- \alpha_n(-Z[X^n]_{\sigma}- 
C[S^n,\beta;X^n]_{\sigma})\}]\} \leq 
1 - \mathbb{E}[\exp\{\hat{\mu} S^{2/(4-\beta)} \cdot \langle Y \rangle_\sigma\}].
\end{equation*}
The upper bound is
attained by the efficient scheme $X^n$ defined
by (\ref{eff}) with $\epsilon_n = \nu \alpha_n^{-1}$ and
\begin{equation*}
\nu =  \mu^{1/2}
\left|\frac{\mu}{6(2-\beta)} \right|^{-1/(4-\beta)}.
\end{equation*}
This can be proved by applying Theorem~2.6 of Fukasawa~\cite{AAP}.
For $\beta \in (1,2)$, similarly  we get
\begin{equation*}
C_\sigma^{2/(2-\beta)}V_{\sigma} \geq \frac{\mu}{18} | S^{2/(4-\beta)}
 \cdot \langle Y \rangle_\sigma |^{(4-\beta)/(2-\beta)}
\end{equation*}
and so,
\begin{equation*}
\begin{split}
C_\sigma + \frac{1}{2}V_{\sigma} &\geq 
C_\sigma + \frac{\mu}{36} | S^{2/(4-\beta)}
 \cdot \langle Y \rangle_\sigma |^{(4-\beta)/(2-\beta)}
C_\sigma^{-2/(2-\beta)}  \\ &\geq
\check{\mu}
S^{2/(4-\beta)}
 \cdot \langle Y \rangle_\sigma,
\end{split}
\end{equation*}
where
\begin{equation*}
\check{\mu} = \left|\frac{\mu}{18(2-\beta)} \right|^{(2-\beta)/(4-\beta)}  +
 \frac{\mu}{36}
\left|\frac{\mu}{18(2-\beta)} \right|^{-2/(4-\beta)}.
\end{equation*}
Therefore,
\begin{equation*}
\lim_{n \to \infty}
\{1- \mathbb{E}[\exp\{- \alpha_n(-Z[X^n]_{\sigma}- 
C[S^n,\beta;X^n]_{\sigma})\}]\} \leq 
1 - \mathbb{E}[\exp\{\check{\mu} S^{2/(4-\beta)} \cdot \langle Y \rangle_\sigma\}].
\end{equation*}
The upper bound is asymptotically
attained by the efficient scheme $X^n = X^n(\gamma)$ defined
by (\ref{gamma}) and (\ref{Xgamma}) when $|\gamma| \to \infty$, where
 $\epsilon_n = \check{\nu} \alpha_n^{-1}$ and
\begin{equation*}
\check{\nu} =  \mu^{1/2}
\left|\frac{\mu}{18(2-\beta)} \right|^{-1/(4-\beta)} \left|
\frac{\cosh((\beta-1)\gamma)}{\cosh(\gamma)}
\right|^{1/(2-\beta)}.
\end{equation*}
Consequently, the 
efficient schemes obtained in the preceding sections
are in fact maximizers of the exponential utility in an asymptotic sense.


\begin{thebibliography}{99}
\bibitem{DK}
Denis, E. and Kabanov, Y. :
Mean square error for the Leland-Lott hedging strategy: convex pay-offs.
{\it Finance Stoch.} 14,  no. 4, 625-667 (2010)

\bibitem{RFE}
Fukasawa, M. :
Asymptotic efficiency for discrete hedging strategies (in Japanese).
{\it Selected papers for the 10 th anniversary of
Financial Technology Research Institute, Inc.} (2009)

\bibitem{AAP}
Fukasawa, M. :
Discretization error of stochastic integrals.
{\it Ann. Appl. Probab.} 21, 1436-1465 (2011)

\bibitem{SAFA}
Fukasawa, M. :
 Asymptotically efficient discrete hedging. 
{\it Stochastic Analysis with Financial Applications}, 
Progress in Probability 65, 331-346 (2011) 
\bibitem{RAFE}
Fukasawa, M. :
Conservative delta hedging under transaction costs.
to appear in {\it Recent Advances in Financial Engineering},  World Scientific (2012) 

\bibitem{GG06}
Geiss, C. and Geiss, S.: 
On an approximation problem for stochastic integrals where random time
	nets do not help.
{\it  Stochastic Process. Appl.} 116, 407-422  (2006)

\bibitem{GT09}
Geiss, S. and Toivola, A.:
Weak convergence of error processes in discretizations of stochastic 
integrals and Besov spaces. {\it Bernoulli} 15, no. 4, 925-954  (2009)

\bibitem{GJ}
Genon-Catalot, V. and Jacod, J.: Estimation of the diffusion coefficient
	for diffusion processes: random sampling. 
{\it Scand. J. Statist.} 21, no. 3, 193-221  (1994)

\bibitem{GT01}
Gobet, E.; Temam, E. :  Discrete time hedging errors for
      options with irregular payoffs. {\it Finance Stoch.}  5, no.3,
      357-367 (2001)

\bibitem{HM2005}
Hayashi, T. and Mykland, P.A. :
Evaluating hedging errors: an asymptotic approach.
{\it Math. Finance} 15, no. 2, 309-343 (2005) 


\bibitem{JS}
Jacod, J. and Shiryaev, A.N.:
{\it Limit theorems for stochastic processes}.
2nd ed., Springer-Verlag  (2002)

\bibitem{Karandikar1995}
Karandikar, R.L. :
 On pathwise stochastic integration.
{\it  Stochastic Process. Appl.}  57,  no. 1, 11-18 (1995)


\bibitem{KS}
Karatzas, I. and Shreve, S.E.:
{\it Brownian Motion and Stochastic Calculus}.
Springer-Verlag,
New York  (1991)


\bibitem{Leland}
Leland, H.E.:
Option pricing and replication with transaction costs.
{\it Journal of Finance} 40, 1283-1301  (1985)



\bibitem{Rootzen}
Rootz\'en, H. :
Limit distributions for the error in approximations of stochastic
        integrals.
 {\it  Ann. Probab.}  8, no. 2, 241-251 (1980)

\bibitem{TV}
Tankov, P. and  Voltchkova, E.:
 Asymptotic analysis of hedging errors in models with jumps. {\it Stochastic
	Process. Appl.} 119, no. 6, 2004-2027 (2009)
\end{thebibliography}
\end{document}